\documentclass[12pt]{amsart}
\usepackage{amsmath}
\usepackage{amscd}
\usepackage{graphics}
\usepackage{latexsym}

\newtheorem{Thm}{Theorem}

\newtheorem{Def}{Definition}

\def\C{{\mathbb C}}

\def\R{{\mathbb R}}
\def\P{{\mathbb P}}

\def\CF{\mathcal F}

\begin{document}

\title[]{ Exotic Structures on smooth $4$-manifolds }

\author{Selman Akbulut }
\thanks{The author is partially supported by NSF grant DMS 9971440 and IMBM}
\address{Department  of Mathematics, Michigan State University,  MI, 48824}
\email{akbulut@math.msu.edu }
\date{\today}
\begin{abstract} 
 A short survey of exotic smooth structutes on $4$-manifolds is given with a special emphasis on the corresponding cork structures. Along the way we discuss some of the more recent results in this direction, obtained jointly with  R. Matveyev, B.Ozbagci, C.Karakurt and  K.Yasui.
 
\end{abstract}
\maketitle

\setcounter{section}{-1}

\vspace{-.3in}

\section{Corks}

Let $M$ be a smooth closed simply connected $4$-manifold, and  $M'$ be an exotic copy of $M$ (a smooth manifold homeomorphic but not diffeomorphic to $M$). Then we can find a  compact contractible  codimension zero submanifold $W\subset M$ with complement $N$, and an involution $f:\partial W\to \partial W$ giving a decomposition:

\begin{equation}
M=N\cup_{id}W \;,  \hspace{0.1in} M'=N\cup_{f}W
\end{equation}

 \begin{figure}[ht]  \begin{center}  
\includegraphics{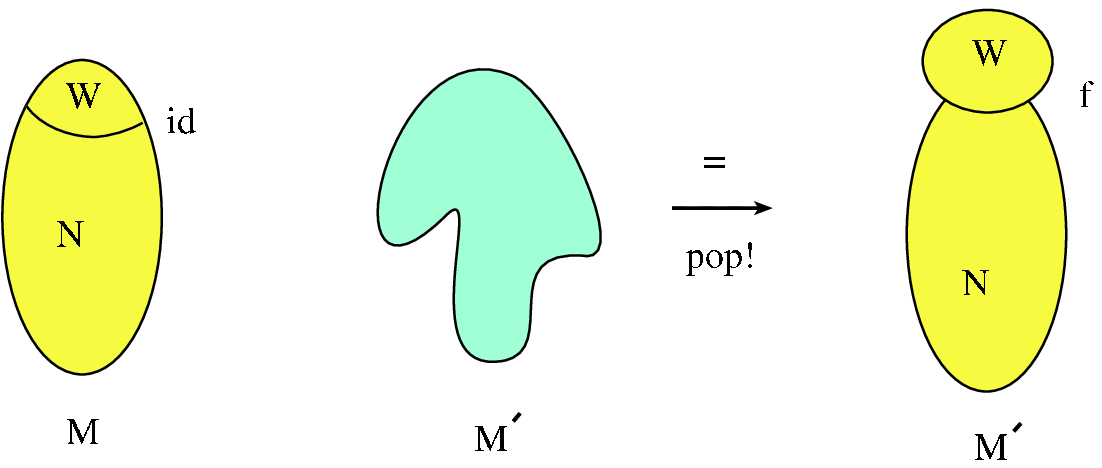}   
\caption{}    \end{center}
  \end{figure}

The existence of this structure was first observed on an example in \cite{a1}, then in  \cite{m} and  \cite{cfhs} it was generalized to the general form discussed above; and another improved version was given in \cite{k}. Since then the contractible pieces $W$ appearing  in this decomposition has come to be known as ``{\it corks}".

\vspace{.05in}

By \cite{am}, in the cork decompositions, each $W$ and $N$ pieces can be made Stein manifolds. 
This is achieved by a useful technique (``creating positrons'' in \cite{am}) which amounts to moving the  common boundary
$\Sigma =\partial W=\partial N $ in $M$ by a convenient homotopy:  
First, by handle exchanges we can assume that each $W$ and $N$ side has only $1$- and $2$-handles. Eliashberg criterion (cf. \cite{g}) says that manifolds with $1$- and $2$-handles are Stein if the attaching framings of the $2$-handles are sufficiently negative (say {\it admissable}), i.e. any $2$-handle $H$ has to be attached along a knot $K$ with framing less than the Thurston-Bennequin framing $tb(K)$ of any Legenderian representation of $K$. The idea is by local handle exchanges near $H$ (but away from $H$) to alter $\Sigma \leadsto \Sigma' $ which results in an increase in Thurston-Bennequin numbers
$tb(K) \leadsto tb(K')=tb(K)+3$.
\begin{figure}[ht]  \begin{center}  
\includegraphics{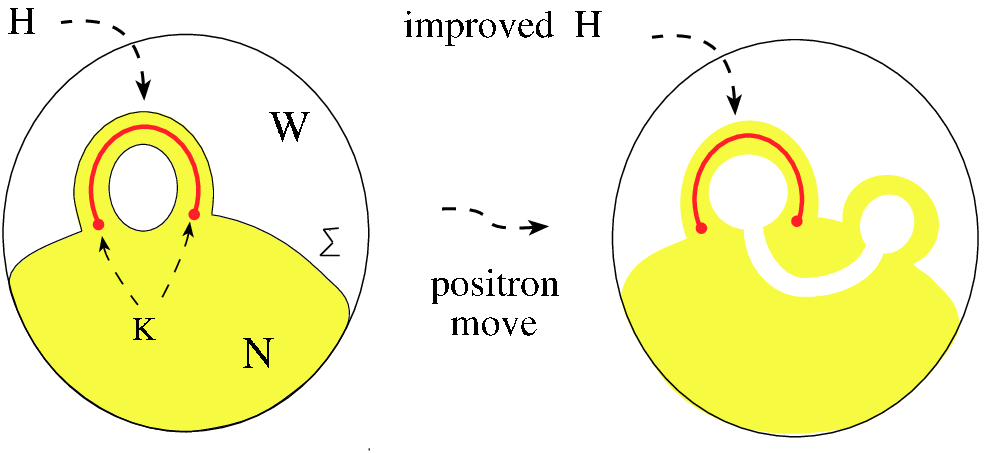}   
\caption{}    \end{center}
  \end{figure}
  
  For example, as indicated in Figure 3, by carving out a tubular neighborhood of a properly imbedded $2$-disc ${\bf a}$  from the interior of $N$ increases $tb(K)$ by $3$.  Carving in the $N$ side corresponds to attaching a $2$-handle ${\bf A}$ from the $W$ side, which itself might be attached with a``bad''  framing. To prevent this, we also attach a $2$-handle ${\bf B}$ to N near ${\bf a}$, which corresponds to carving out a $2$-disk from the $W$ side. This makes the framing of the $2$-handle ${\bf A}$ in the $W$ side admissable, also ${\bf B}$ itself is admisable. So by carving a $2$-disc ${\bf a}$ and attaching a $2$-handle ${\bf B}$ we improved the attaching framing of the $2$-handle $H$, without changing other handles (we changed $\Sigma$ by a homotopy). This technique gives: 

    \begin{figure}[ht]  \begin{center}  
\includegraphics{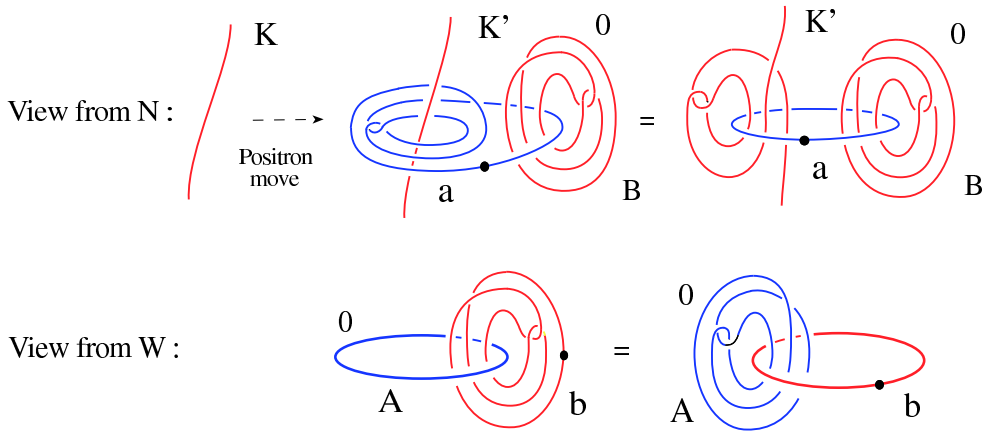}   
\caption{Making W and N Stein}    \end{center}
  \end{figure}

\vspace{.03in}

{\Thm (\cite {am}) Given  any  decomposition  
 of a  closed smooth $4$-manifold $M=N_1\cup _{\partial} N_2$ by codimension zero submanifolds, with each piece consisting of $1$-and $2$-handles; after altering pieces by a homotopy we can get a similar decomposition $ M=N_1'\cup N_2'$ where both pieces $N_1$ and $N_2$ are Stein manifolds.}

\vspace{.04in}
By applying \cite{gi} one can also assume that the two open books on the common contact $3$-manifold boundaries  match, but with the wrong orientation (\cite{b}).

 {\Def  A  Cork is a pair $(W,f)$, where $W$ is a compact Stein manifold, and $f:\partial W\to \partial W$ is an involution, which extends to a self-homemorphism of $W$, but it does not extend to a self-diffeomorphism of $W$. We say $(W,f)$ is a cork of $M$, if we have the decomposition (1) for some exotic copy $M'$ of $M$.}
 
 \vspace{.05 in} 
    
 In particular a cork is a fake copy of itself. There are some natural families of corks $W_n$, $n=1,2,..$  which are generalization of the Mazur manifold $W$ used in \cite{a1}, and $\bar{W}_n$, $n=1,2,..$ which were introduced in \cite{am} (so called positrons). It is a natural question to ask  whether these small standard corks are sufficient to explain all exotic smooth structures on $4$-manifolds ? (Figure 5).  For example, $W_1$ is a cork of $E(2)\#\bar{\C\P^2}$ (\cite{a1}), and $\bar{W_{1}}$ is a cork of the Dolgachev surface $E(1)_{2,3}$ (\cite{a2}), where $E(n)$ is the Elliptic surface of signature $-8n$.
 
  \begin{figure}[ht]  \begin{center}  
\includegraphics{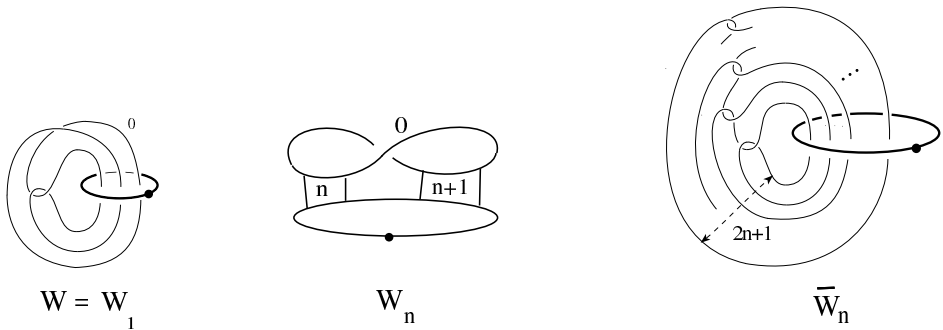}   
\caption{Variety of  corks}    \end{center}
  \end{figure}

   \begin{figure}[ht]  \begin{center}  
\includegraphics{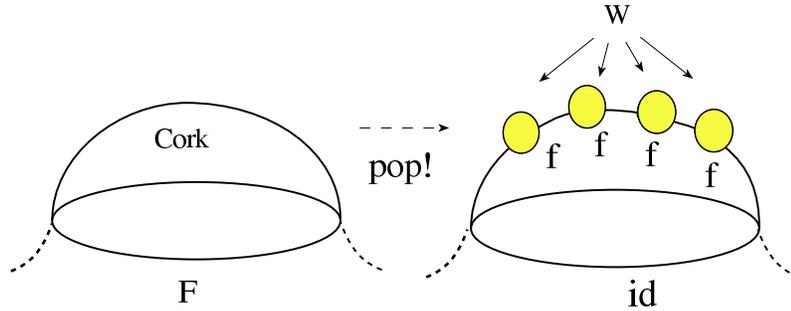}   
\caption{Do all corks decompose into standard corks?}    \end{center}
  \end{figure}

  \vspace{.05 in} 
  
 The reason we require corks to be Stein manifolds is to rule out trivial examples, as well as introduce rigidity in their structures. For example, a theorem of Eliashberg says that if a Stein manifold has boundary $S^3$ or $S^1\times S^2$ then it has to be a $B^4$ or $S^1\times B^3$, respectively.

\subsection{ How to recognize a cork?}

$\:$

\vspace{.1in}

In general it is hard to recognize when a codimension zero contractible submanifold $W\subset M^4$  is a cork of $M$. In fact, all the  corks obtained in the general cork decomposition theorem of \cite{m} has the property that $W\cup -W=S^4$ and $W\cup_{f} -W=S^4$. So it is easy to imbed $W$'s into charts of $M$ without being corks of $M$. One quick way of showing $(W,f)$ is a cork of a manifold $M$ with nontrivial Seiberg-Witten invariants, is to show that the change  $M\leadsto M'$ in (1) gives a split manifold $M'$, implying zero (or different) Seiberg-Witten invariants. In  \cite{ay1} and \cite{a2}, by using this strategy many interesting corks were located. 

\vspace{.05in} 

There are also some hard to calculate algebraic ways of checking if $W\subset M $ is a cork, provided we know the Heegard-Floer homology groups of the boundary of $W$ \cite{os}.   This follows from the computation of the Ozsvath-Szabo $4$-manifold invariant, i.e. by first removing two $B^4$'s from $M$ as shown in Figure 6, and computing certain trace of the induced map on the Floer homology of the two $S^3$ bounday components  (induced from the cobordism). For example we have:

 \begin{figure}[ht]  \begin{center}  
\includegraphics{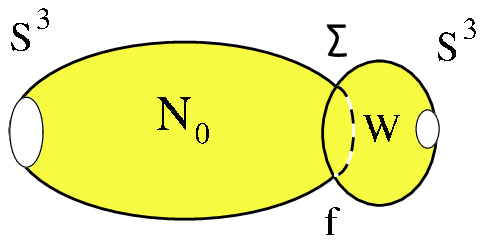}   
\caption{}    \end{center}
  \end{figure}

{\Thm (\cite{ad}): Let $M=N\cup _{\partial} W$ be a cork decomposition of smooth closed $4$-manifold, where $W$ is the Mazur manifold, and $b_{2}^{+}(M)>1$ (union is along the common boundary $\Sigma$). Let $N_{0}$, be the cobordism from $S^3$ and $\Sigma$  obtained from $N$ by removing a $B^4$ from its interior. Then $Q'= N\cup_{f} W$ is a fake copy of $Q$, if the image of the following map lies in $T^{+}_{0}$ for some $Spin^c$ structure $s$.
$$F^{mix}_{(N_{0},s)}: HF^{-} (S^3) \to HF^{+}(\Sigma)\cong T_0^+ \oplus {\bf Z}_{(0)} \oplus {\bf Z}_{(0)}$$.}

 \subsection{Constructing exotic manifolds from Corks}

$\:$

\vspace{.1in}

By thickening a corks in  two different ways one can obtain absolutely exotic manifolds pairs (i.e. homeomorphic but not diffeomorphic manifolds). Here is a quick review of \cite{a3}:  Let  $(W,f)$ be the Mazur cork, the  $f:\partial W \to \partial W$ has an amazing property: There are pair of loops $\alpha $, $\beta$ with:

\begin{itemize}
\item $f(\alpha)=\beta$

\item $M := W + (2\mbox{\it - handle to}\; \alpha \; \mbox{\it with} -1 \;\mbox{framing})  \;$  is a Stein manifold.
\item $\beta $ is slice in $W$, hence \\
$M' := W + (2\mbox{-handle to}\; \beta \; \mbox{\it with} -1 \;\mbox{framing})$ \\ contains an imbedded $-1$ sphere 
\end{itemize}

So $M'$ is an obsolutely exotic copy of $M$; if not $M'$ would be a Stein manifold also, but any Stein manifold compactifies into irreducible symplectic manifold (\cite{lm}), contradicting the existence of the smoothly imbedded $-1$ sphere.

  \begin{figure}[ht]  \begin{center}  
\includegraphics{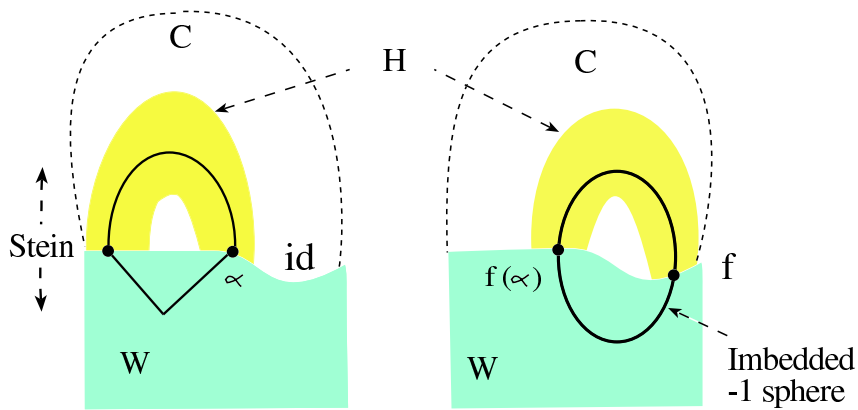}   
\caption{Inflating a cork to exotic manifold pairs}    \end{center}
  \end{figure}

Interestingly, by handle slides one can show that each of  $M$ and $M'$ is obtained by attaching a $2$- handle to $B^4$ along a knot, as indicated in Figure 8  \cite{a3}. 
  \begin{figure}[ht]  \begin{center}  
\includegraphics{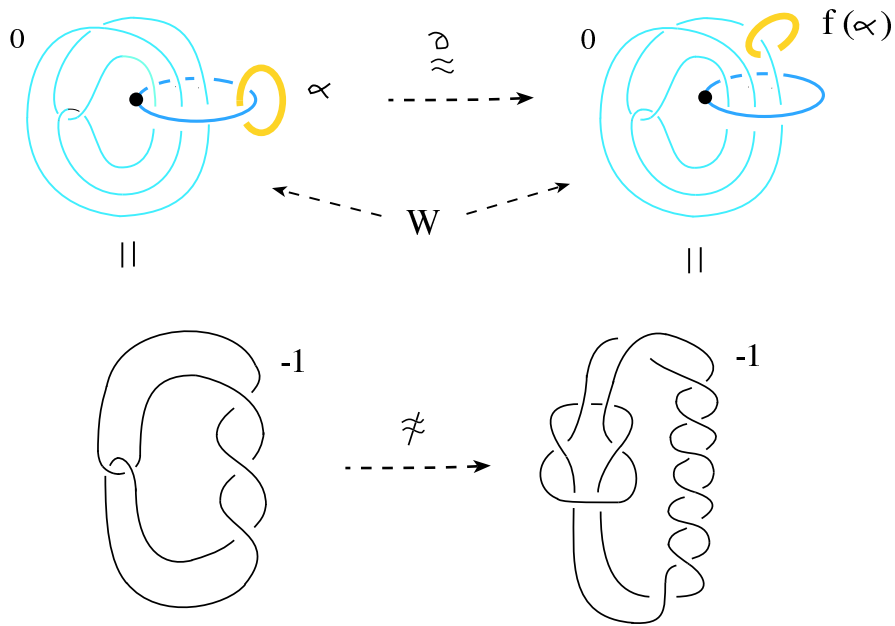}   
\caption{Exotic manifold pairs}    \end{center}
  \end{figure}

Reader should contrast this with  \cite{a5}, where examples of other knot pairs $K,L\subset S^3$ are given (one is slice the other one is not slice) so that attaching $2$-handles  to $B^4$ gives diffeomorphic $4$-manifolds, as in Figure 9.
  \begin{figure}[ht]  \begin{center}  
\includegraphics{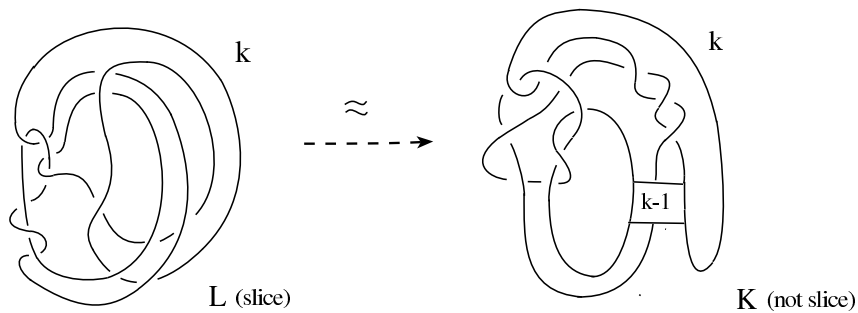}   
\caption{Diffeomorphic manifold pairs}    \end{center}
  \end{figure}

\section{PALFs} 

It turns out that Stein manifolds admit finer structures as  primitives,  we called them PALF's in \cite{ao1}  $ Stein = |PALF|$: ``Positive Allowable Lefschetz Fibration'' over the $2$- disk, where the regular fibers are surfaces $F$ with boundaries. Here``allowable'' means that the monodromies of Lefschetz singularities over the singular points are products of positive Dehn twists along non-separating loops (this last condition is not a restriction, it comes for free from the proofs). Existence of this structure on Stein manifolds was first proven in \cite{lp}, later in \cite{ao1} a constructive topological proof along with its converse is given, hence establishing:

{\Thm   There is a surjection
$$ \left\{ PALF's \right\} \Longrightarrow \left\{ \mbox{ Stein Manifolds} \right\} $$}
Here by Eliashberg's characterization (\cite{g}), a {\it Stein Manifold} means a handlebody consisting of $1$- and $2$-handles, where the $2$-handles are attached along Legenderian framed link, with each component K framed with  $ tb(K)-1$  framing.
 
  \begin{figure}[ht]  \begin{center}  
\includegraphics{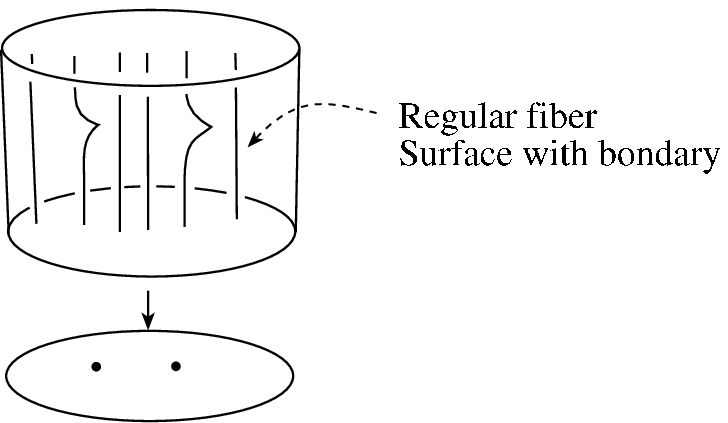}   
\caption{PALF}    \end{center}
  \end{figure}

The proof that a PALF $\CF$  gives a Stein manifold $ |\CF| $ goes as follows: By \cite{ka} $\CF$ is obtained by starting with the trivial fibration $X_0=F\times B^2 \to B^2$ and attaching sequence of $2$-handles to the curves $k_i\subset F$ , $i=1,2..$ on the fibers, with framing  one less that the  page framing: 
$X_0 \leadsto X_1\leadsto .. \leadsto X_n=\CF$. On a $F \times B^2$ we start with the standard Stein structure and assume on the contact boundary $F$ is a convex surface with the ``divider'' $\partial F$ \cite{t}, then by appying the  ``Legenderian realization principle"  of \cite{h}, after an isotopy we make the surface framings of  $k_i\subset F$ to be the Thurston-Bennequin framings, and then the result follows from Eliashberg's theorem. 

\vspace{.05in}

Conversely, to show that a Stein manifold $W$ admits a PALF,  (here we only indicate the proof when there are no $1$-handles) we isotope the Legenderian framed link to square bridge position (by turning each component counterclockwise 45 degrees), and put the framed link on a fiber $F$  of the $(p,q)$ torus knot $L$ as indicated  in Figure 11. $L$ gives a PALF structure on $B^4=|\CF|$. Now attaching handles to this framed link has the affect of enhancing the monodromy of the $(p,q)$ torus knot by the Dehn twist along them, resulting in a bigger PALF.  An improved version of this theorem is given in \cite{ar}

  \begin{figure}[ht]  \begin{center}  
\includegraphics{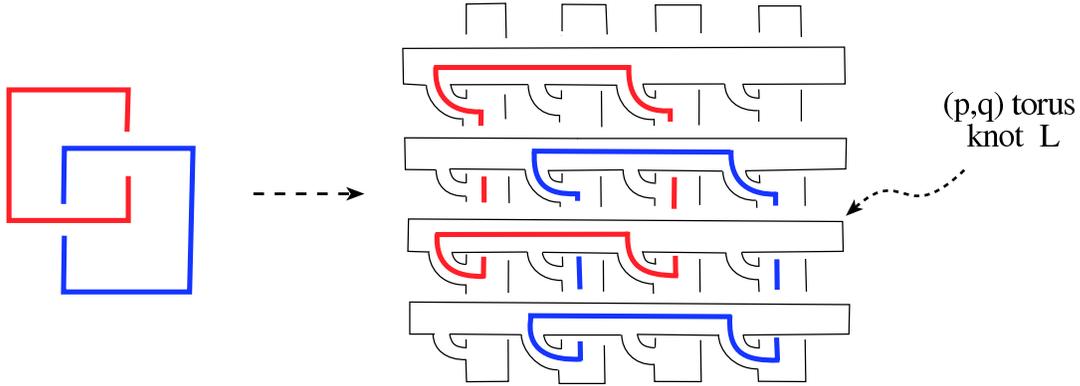}   
\caption{Surgering framed link induces Dehn twists on the page}    \end{center}
  \end{figure}

A PALF structure  $X=|\CF|$ should be viewed as an  auxiliary topological structures $\CF$ on a Stein manifold $X$, like a triangulation or handlebody structures on a smooth manifold. On the boundary  a PALF gives an open book compatible with the induced contact manifold $\partial X=|\partial \CF |$. Usually  geometric  structures come as primitives of  topological structures: $\mbox{Topology}=|\mbox{Geometry}|$, like the real algebraic structures or complex structures on a smooth manifold; but surprisingly in this case the roles are reversed:  $ \mbox{Geometry} =|\mbox{Topology}| $. For example $B^4$ has a unique Stein structure, whereas it has infinitely many PALF structures corresponding to fibered links.

\vspace{.05in}

Choosing an underlying PALF  is often useful in solving  problems in Stein manifolds.  A striking application of this principle was the approach in \cite{ao2} to the compactification problem of Stein manifolds, which was later strengthen by \cite{e} and  \cite{et}. The problem of  compactifying a Stein manifold $W$ into a closed symplectic manifold was first solved in \cite{lm}, then in \cite{ao2} by using PALF's  an algorithmic solution was given. An analogous case is compactifying  into closed manifold; by first choosing a handlebody on $W$, then canonically compactifying it by attaching dual handles (doubling). In the symplectic case we first choose a PALF on $W=|\CF|$, then attach a $2$-handle to the binding of the open book on the boundary (Figure 12) and get a closed surface $F^{*}$-bundle over the $2$-disk with monodromy product of positive Dehn twists $\alpha_1.\alpha_2.. \alpha_k$. We then extend this fibration $\CF$ by doubling monodromies $\alpha_1,\alpha_2..\alpha_k \alpha_k^{-1}..\alpha_{1}^{-1}$ (i.e. attaching corresponding $2$-handles) and capping off with $F^{*}\times B^2$ on the other side. We do this after converting each negative Dehn twist $\alpha^{-1}_{i}$ in this expression into products of positive Dehn twists by using the relation $(a_1b_1...  a_gb_g)^{4g+1}=1$ among the standard Dehn twist generators of the  surface $F^{*}$ of genus $g$ (cf.\cite{ao2}). 
 
 \begin{figure}[ht]  \begin{center}  
\includegraphics{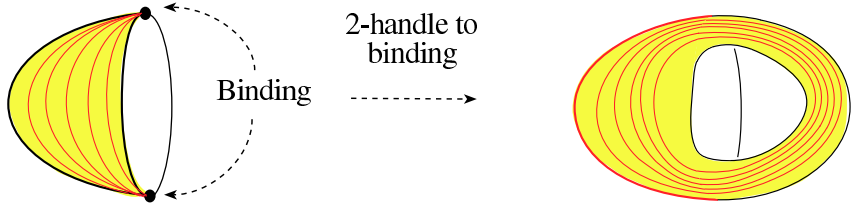}   
\caption{Adding a $2$-handle to binding}    \end{center}
  \end{figure}

  \vspace{-.3in}
  
     \begin{figure}[ht]  \begin{center}  
\includegraphics{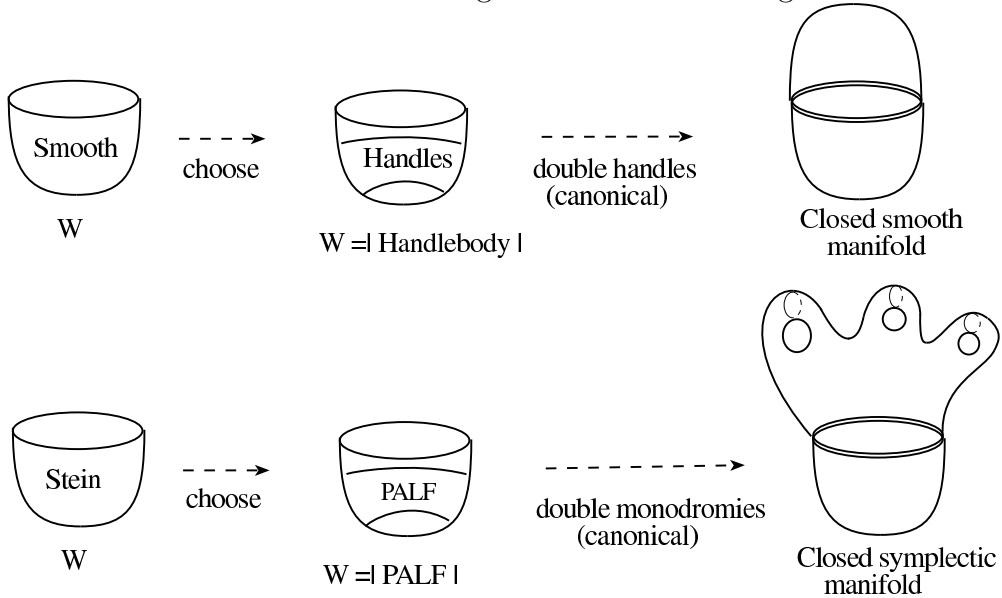}   
\caption{Natural compactification after choosing an auxiliary  structure}    \end{center}
  \end{figure}

\noindent After choosing a PALF on $W$, this becomes an algorithmic canonical process, even in the case of Stein ball $B^4$ by using different PALF structures on it we get variety of different symplectic compactifications of $B^4$, for example: 

 \vspace{.05in}

 $ B^4 =| \mbox{Unknot} | \leadsto S^2 \times S^2 $, whereas
 $\; \;B^4 =| \mbox{Trefoil} | \leadsto  K3$ Surface.

\section{BLFs}

Lack of a uniqueness result in the ``cork decompositions'', and the differing orientations of the two Stein pieces obtained from  Theorem 1 makes it hard to define $4$-manifold invariants.  In \cite{adk} a more general version of the Lefschetz fibration (or pencil) structure on $4$-manifolds introduced, namely  ``Broken Lefschetz fibration"  BLF (or  Broken Lefschetz Pencil ``BLP''), where they allow Lefschetz fibrations (or pencils) $\pi:  X^4\to S^2$ to have circle singularities, i.e. on a neighborhood of some circles   $\pi$ can look like a map $S^1\times B^3\to \R^2$ given by $(t,x_1,x_2,x_3)\mapsto (t, x_1^2+x_2^2-x_3^2)$ (otherwise it is a Lefschetz fibration). In \cite{adk}, by using analytic techniques, it was shown that every $4$-manifold $X$ with $b_{+}^{2}>0$ is a BLP. Also, after a useful partial result in \cite{gk}, in \cite{l} and  \cite{ak} two independent proofs that all $4$-manifolds are BLF have been given; the first proof uses singularity theory and the second uses the handlebody theory. Now in \cite{p} there is an approach to use BLF's to construct $4$-manifold invariants.

\vspace{.05in}

The proof of \cite{ak} proceeds along the lines of Theorems 1 and 3, discussed earlier, roughly it goes as follows: First define ALF's which are weaker version of PALF's, they are ``Achiral Lefschetz Fibrations" over the $2$-disk with bounded fibers, where we allow Lefschets singularites to have monodromies that are negative Dehn twists. From the proof of  \cite{ao1} we get a surjection:
$$ \left\{ ALF's \right\} \Longrightarrow \left\{ \mbox{ Almost Stein Manifolds} \right\} $$
Here an {\it almost Stein Manifold} means a handlebody consisting of $1$- and $2$-handles, where the $2$-handles are attached to a Legenderian framed link, with each component $K$ framed  with $tb(K)\pm 1$  framing (it turns out every $4$ dimensional handlebody consisting of $1$-and $2$-handles has this nice structure). 

\vspace{.05in}

First, we make a tubular neighborhood $X_2$ of any imbedded surface in $X$ a ``concave BLF'' (e.g. \cite{gk}). A {\it concave BLF} means a BLF with $2$- handles attached to circles transversal to the pages on the boundary (open book), as indicated in Figure 15 (so a concave BLF fibers over the whole $S^2$, with closed surface regular fibers on one hemisphere, and $2$-disks regular fibers on the other hemisphere).  Also we make sure that the complement $X_1=X-X_{2}$ has only $1$-and $2$-handles, hence it is an ALF.  By applying \cite{gi} we make sure that the boundary open books induced from each side $X_i$, $i=1,2$ match. So we have a (matching) union $X=X_1 \cup X_2$, consisting of an ALF $X_1$ and concave BLF $X_2$. This would have made X a BLF if $X_1$ were a PALF. 

\vspace{.05in}

Now comes the crucial point.   Like getting a butterfly from a silkworm by drilling its cocoon, we will turn the ALF  $X_1$ into a PALF by removing a disk from it, i.e. we will apply the positron move of  Theorem 1: For each  framed knot representing  a $2$-handle of $X_1$, we pick an unknot  ($K$ in Figure 14) with the properties: (1)  It lies on a page, (2) It links  that framed knot twice, as in Figure 14. These conditions allows us to isotope K to the boundary of a properly imbedded $2$-disk in $X_1$ meeting each fiber of the ALF once (here K is isotoped to the meridian of the binding curve). Therefore carving out  the tubular neighborhood of this disk from $X_1$  preserves the ALF structure (this is indicated by putting a dot in $K$ in the figure, a notation from \cite{a5}), where the new fibers are obtained by puncturing the old ones. But now, this magically changes the framings of each framed knots by $tb(K)+1 \leadsto tb(K)-1$, hence after carving out $X_1$ we end up getting a PALF.  On the other $X_2$ side, this carving corresponds to attaching a $2$ handles to $X_2$ to circles transverse to fibers, so the enlarged $X_2$ is still a concave BLF, so the open books of each side still match.

 \vspace{.1in}
 
  \begin{figure}[ht]  \begin{center}  
\includegraphics{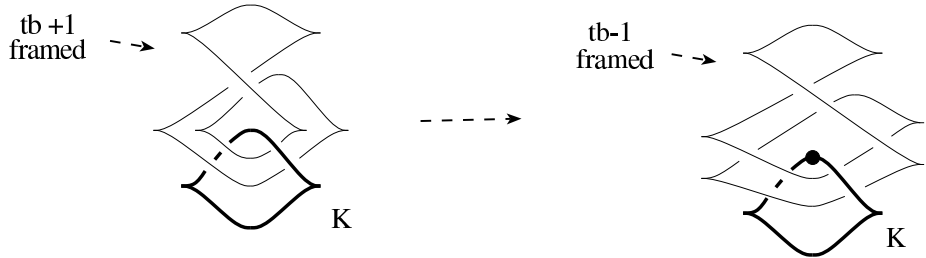}   
\caption{Turning ALF to PALF by carving}    \end{center}
  \end{figure}

  \begin{figure}[ht]  \begin{center}  
\includegraphics{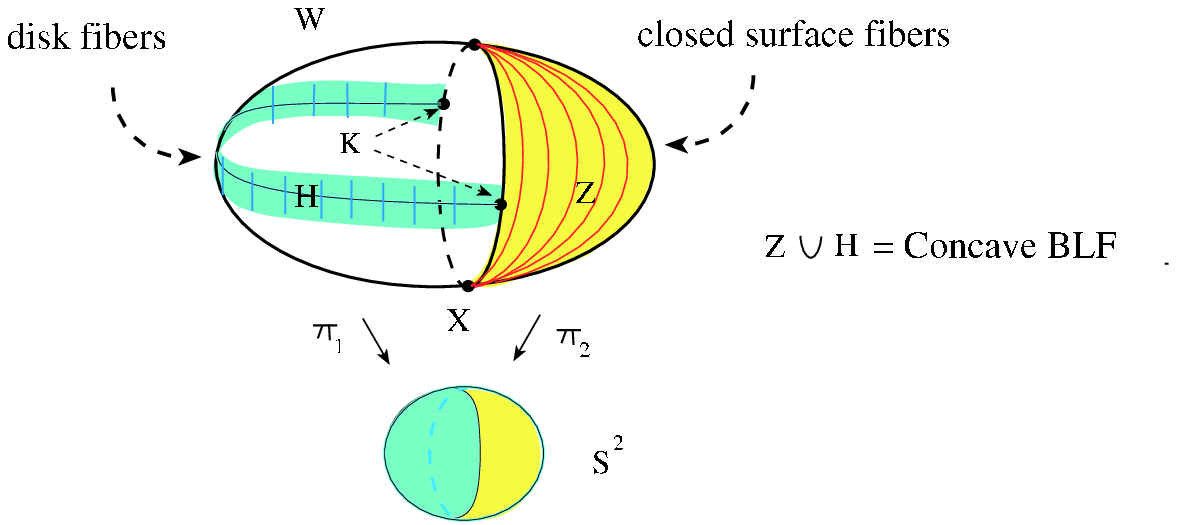}   
\caption{}    \end{center}
  \end{figure}

\newpage

\section{Plugs}

To understand exotic structures of $4$-manifolds better, recently  in \cite{ay1} Yasui and I started to search for   corks in the known examples of exotic $4$-manifolds, since we know theoretically that they exist (Section 0). From this endeavor  we learned two important lessons: $1$- Similar to corks there are different behaving codimension zero submanifolds which are also responsible for the exoticness of $4$-manifolds (we named them {\it plugs}), $2$- The position of corks and plugs in $4$-manifolds plays an important role, for example it helps us to construct exotic Stein manifolds in Theorem 4, whose existence have eluded us for a long time. In some sense plugs generalize the {\it Gluck twisting} operation, just as corks generalize the Mazur manifold. The rest of this section is a brief summary of \cite{ay1}, \cite{ay2}.
 
 {\Def A Plug is a pair (W,f), where $W$ is a compact  Stein manifold, and $f:\partial W\to \partial W$ is an involution, which does {\bf not} extend to a self-homemorphism of 
 $W$, and there is the decomposition (1) for some exotic copy $M'$ of $M$.}

\vspace{.1in}

 Plugs might be deformations of corks (to deform corks to each other we might have to go through plugs). We can think of corks and plugs as freely moving particles in $4$-manifolds like``Fermions" and ``Bosons" in physics (Figure 21), or little knobs on a wall to turn on and off the ambient exotic lights in a room.

\vspace{.05in}

An example of a plug which frequently appears in $4$-manifolds is $W_{m,n}$, where $m\geq1, n\geq 2$ (Figure 16). By canceling the $1$-handle and the $-m$ framed $2$-handle, we see that $W_{m,n}$ is obtained from $B^4$ by attaching a $2$-handle to a knot with $-2n-n^2m^2$ framing. The involution $f$ is induced from the symmetric link.

  \begin{figure}[ht]  \begin{center}  
\includegraphics{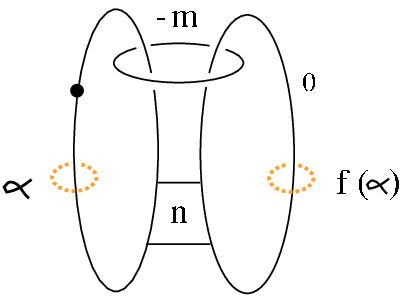}   
\caption{$W_{m,n}$}    \end{center}
  \end{figure}

Here the degenerate case is also interesting. By removing $W_{1,0}$ and gluing with $f$  corresponds to the {\it Gluck twisting operation}. Up to now we only know one example of an exotic manifold which is obtained from the standard one by the Gluck operation, and that manifold is non-orientable  \cite{a4}.

\vspace{.05in}

Notice, if this obvious involution  $f:\partial W_{m.n} \to \partial W_{m,n}$  extended to a homeomorphism, we would get homeomorphic manifolds $W_{m,n}^{1}$ and $W_{m.n}^{2}$, obtained by attaching  
$2$-handles to $\alpha$ and $f(\alpha)$ with $-1$ framings, respectively. But $W_{m,n}^{2}$ and $W_{m,n}^{1}$ have the following non-isomorphic intersection forms, contradiction.
 \begin{equation*}
\left(
\begin{array}{cc}
-2n-mn^2 &1  \\      
1 &-1 
\end{array}
\right), 
\left(
\begin{array}{cc}
-2n-mn^2 &-1-mn  \\
-1-mn &-1-m 
\end{array}
\right), 
\end{equation*}

\vspace{.1in}

The following theorem implies that  $(W_{m,n},f)$ is a plug, and also it says that this plug can be inflated to exotic Stein manifold pairs.

{\Thm (\cite{ay2}) The following simply connected Stein manifolds are exotic copies of each other. }

   \begin{figure}[ht]  \begin{center}  
\includegraphics{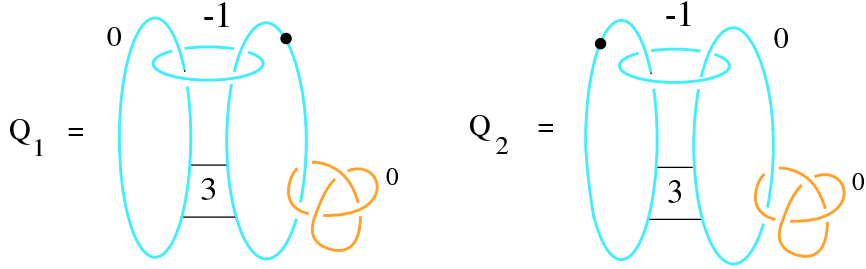}   
\caption{Inflating a plug to an exotic Stein manifold pair}    \end{center}
  \end{figure}
  
 Notice that the transformation  $Q_1 \leadsto Q_2$ is obtained by twisting along the plug $(W_{1,3},f)$ inside. It is easy to check that both are Stein manifolds, and clearly the boundaries of $Q_1$ and $Q_2$ are diffeomorphic, and the boundary diffeomorphism extends to a homotopy equivalence inside, since they have isomorphic  intersection forms $(-1)\oplus (1)$ cf. \cite{bo}. Hence by the Freedman's theorem they are homeomorphic.  The fact that they are are not diffeomorphic follows from an interesting imbedding $Q_1\subset E(2)\#2\bar{\C\P^{2}}$ (so position of plugs are important!), where the two homology generators $<x_1,x_2>$ of $H_2(Q_1)$ intersect the Basic class $K=\pm e_1 \pm e_2$ of  $E(2)\#2\bar{\C\P^{2}}$  with $x_i.e_j=\delta_{ij}$, (here $e_j$ are the two $\bar{\C\P^1}$ factors). Now by applying the adjunction inequality we see that there is no imbedded torus of self intersection zero in $Q_1$; whereas there is one in $Q_2$.

\vspace{.1in}

We can also inflate corks to exotic Stein manifold pairs (compare these examples to the construction in Section 0.2). 

\vspace{.05in}

{\Thm (\cite{ay2}) The following simply connected Stein manifolds are exotic copies of each other. } 

  \begin{figure}[ht]  \begin{center}  
\includegraphics{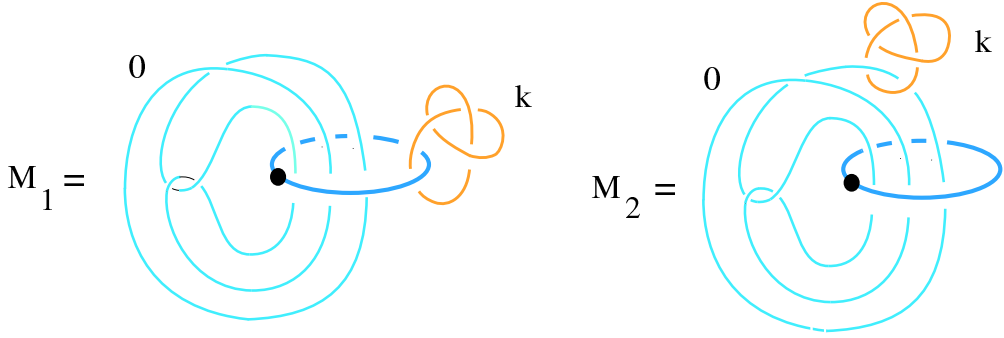}   
\caption{Inflating a cork to an exotic Stein manifold pair ($k\leq 0)$}    \end{center}
  \end{figure}

The proof of this is similar to the previous theorem, the crucial point is finding an imbedding of $M_1$ into a useful closed $4$-manifold with nontrivial Seiberg-Witten invariant.  Note that by a result of Eliashberg, the only Stein filling of $S^3$ is $B^4$.  The existence of simply connected exotic Stein manifold pairs have been first established recently 
in \cite{aems} by using the technique of ``knot surgery''.
Now a natural question: Are all exotic structures on $4$-manifolds induced from the corks $W_n$, $\bar{W_{n}}$ and  plugs $W_{m,n}$? Here is a result of some  checks:

 {\Thm (\cite{ay1})  For $k,r\geq 1$, $n, p,q\geq 2$ and $gcd (p,q)=1$ we have:
 
\vspace{.05 in}
\begin{itemize} 

 \item $E(2k)\# \bar{\bf CP}^2$ has corks  $(W_{2k-1},f_{2k-1})$ and  $(W_{2k},f_{2k})$\\
 
\item $E(2k)\# \; r \bar{\bf CP}^2$ has plugs  $(W_{r, 2k},f_{r,2k})$
and  $(W_{r, 2k+1},f_{r,2k+1})$\\

\item $E(n)_{p,q} \# \bar{\bf CP}^2$ has  cork  $W_{1}$, and has  plug  $W_{1,3}$\\
    
 \item $E(n)_K \# \bar{\bf CP}^2$  has  cork $W_1$, and has  plug $W_{1,3}$. \\
 
\item Yasui's exotic $E(1)\#  \bar{\bf CP}^2$ in \cite{y} has  cork $W_1$.

\end{itemize} }
\vspace{.1 in}

 An interesting question is: Are any two cork imbeddings $(W,f)\subset M$  isotopic to each other?  Put another way, can you knot corks inside of $4$-manifolds?. It turns out that the answer is yes, there are knotted corks  \cite{ay1}. For example, there are  two non-isotopic cork imbeddings $(W_{4},f )\subset M=\C\P^2\# 14 \bar{\C\P^2}$. This is proven by calculating the change in Seiberg-Witten invariants of the two manifolds obtained by twisting $M$ along the two imdedded corks and getting different values, this calculation uses the techniques of \cite{y}.

   \begin{figure}[ht]  \begin{center}  
\includegraphics{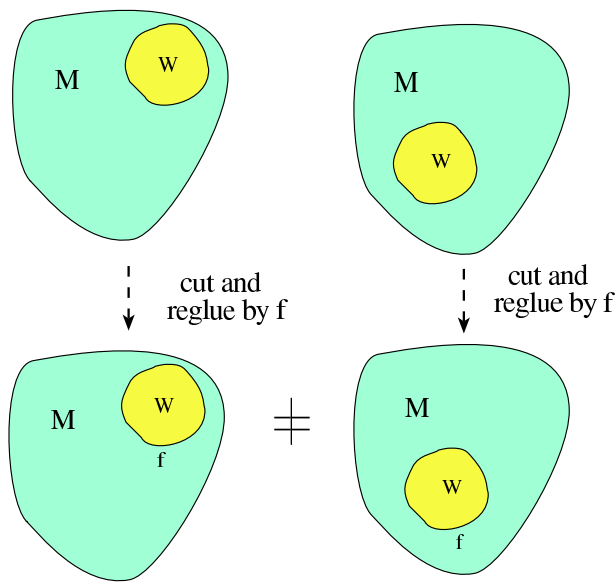}   
\caption{Nonisotopic corks}    \end{center}
  \end{figure}

Another natural question arises: Since every exotic copy of a closed $4$-manifold can be explained by a cork twisting, and there are many ways of constructing exotic copies of  $4$-manifolds, e.g. logarithmic transform, rational blowing down, knot surgery operations (\cite{fs1},\cite{fs2}, \cite{gs}),  are there ways of linking all these constructions to corks? In some cases this can be done  for the rationally blowing down operation $X\leadsto X_{(p)}$, by showing $X_{(p)}\#(p-1)\bar{\bf CP}^2$ is obtained from $X$ by a cork twisting  along some $W_{n}\subset X$ \cite{ay1}. It is already known that there is a similar relation between logarithmic transforms and the rational blowing downs. The difficult remaining case seems to be the problem of relating a general knot surgery operation to cork twisting.

\vspace{.1in}

\noindent {\bf Remark}: {\it We don't yet know if an exotic copy of a manifold with boundary differs from its standard copy by a cork.  It is likely that there is some relative version of the cork theorem. Perhaps the most interesting example to check is the exotic Cusp of \cite{a6}, which is the smallest example of an exotic smooth simply connected manifold we know, which requires $1$- or $3$-handles in any of its handlebody decomposition, whereas its standard copy has only $2$-handles \cite{ay2}. }

 \begin{figure}[ht]  \begin{center}  
\includegraphics{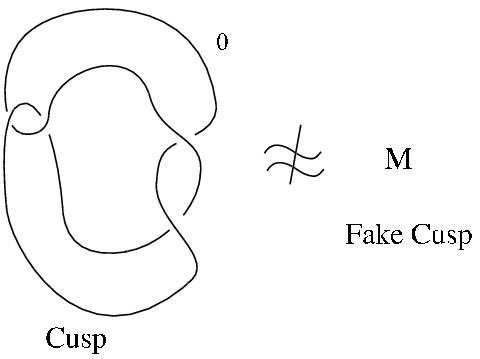}   
\caption{}    \end{center}
  \end{figure}

  \begin{figure}[ht]  \begin{center}  
\includegraphics{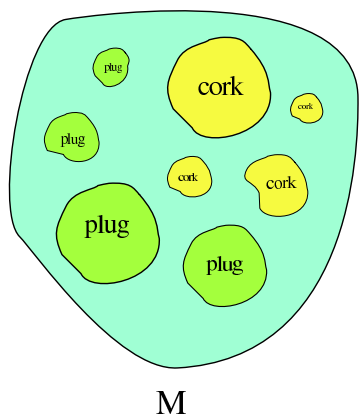}   
\caption{Zoo of corks and plugs in a $4$-manifold}    \end{center}
  \end{figure}

\end{document}